\documentclass[12pt,reqno]{amsart}

\usepackage[arrow,matrix,curve]{xy}

\usepackage[dvips]{graphicx} 

\usepackage{amssymb, latexsym, amsmath, amscd, array, hyperref
}

\usepackage{breakurl}   

\usepackage{enumitem} 

\theoremstyle{definition}

\numberwithin{equation}{section}

\DeclareMathOperator{\adequal}{\;\raisebox{-3pt}{$\ulcorner\!\urcorner$}\;}

\numberwithin{equation}{section}

\author{Mikhail G. Katz}\address{M. Katz, Department of Mathematics,
Bar Ilan University, Ramat Gan 52900
Israel}\email{katzmik@macs.biu.ac.il}

\title[Mathematical conqueror, Unguru polarity, task of history]
{Mathematical conquerors, Unguru polarity, and the task of history}

\begin{document}

\thispagestyle{empty}


\begin{abstract}
We compare several approaches to the history of mathematics recently
proposed by Bl{\aa}sj\"o, Fraser--Schroter, Fried, and others.  We
argue that tools from both mathematics and history are essential for a
meaningful history of the discipline.  

In an extension of the Unguru--Weil controversy over the concept of
\emph{geometric algebra}, Michael Fried presents a case against both
Andr\'e Weil the ``privileged observer'' and Pierre de Fermat the
``mathematical conqueror.''  We analyze Fried's version of Unguru's
alleged polarity between a historian's and a mathematician's history.
We identify some axioms of Friedian historiographic ideology, and
propose a thought experiment to gauge its pertinence.  

Unguru and his disciples Corry, Fried, and Rowe have described
Freudenthal, van der Waerden, and Weil as Platonists but provided no
evidence; we provide evidence to the contrary.  We analyze how the
various historiographic approaches play themselves out in the study of
the pioneers of mathematical analysis including Fermat, Leibniz,
Euler, and Cauchy.
\end{abstract}

\maketitle
\tableofcontents

\section{Introduction}
\label{s1}

The recent literature features several approaches to the history of
mathematics.  Thus, Michael N. Fried (\cite{Fr18}, 2018) and
Guicciardini (\cite{Gu18}, 2018) argue for versions of Unguru's
approach (see below).  Bl{\aa}sj\"o (\cite{Bl14}, 2014) advocates a
rational history as opposed to an ``idiosyncraticist'' one.  Fraser
and Schroter propose something of a middle course that defines the
task of the history of mathematics as ``our attempt to explain why a
certain mathematical development happened'' in (\cite{Fr19a}, 2019,
p.\;16).  We illustrate the latter approach in Section~\ref{s43}, in
the context of certain developments in mathematical analysis from
Euler to Cauchy, following Fraser--Schroter (\cite{Fr20}, 2020).

We analyze the perception of mathematical historiography that posits a
polarity between a historical and a mathematical view.  Such a
perception is often associated with Sabetai Unguru.  Against such
Unguru polarity, we argue that tools from both disciplines are both
useful and essential.

\subsection{Unguru, Weil, van der Waerden, Freudenthal}

Sabetai Unguru (\cite{Un75}, 1975) and Andr\'e Weil (\cite{We78},
1978) famously battled one another over the relation between Greek
mathematics and the concept of \emph{geometric algebra}, a term
introduced by H. G. Zeuthen in 1885 (see Bl{\aa}sj\"o \cite{Bl16},
2016, p.\;326; H{\o}yrup \cite{Ho16}, 2016, pp.\;4--6).  

We note that B.\;L.\;van der Waerden (\cite{Va76}, 1976) and Hans
Freudenthal (\cite{Fr77}, 1977) published responses to Unguru earlier
than Weil.  A clarification is in order concerning the meaning of the
term \emph{geometric algebra}.  Van der Waerden explained the term as
follows:
\begin{quote}
We studied the wording of [Euclid's] theorems and tried to reconstruct
the original ideas of the author.  We found it \emph{evident} that
these theorems did not arise out of geometrical problems.  We were not
able to find any interesting geometrical problem that would give rise
to theorems like II 1--4.  On the other hand, we found that the
explanation of these theorems as arising from algebra worked well.
Therefore we adopted the latter explanation.  Now it turns out
{\ldots} that what we, working mathematicians, found evident, is not
evident to Unguru.  Therefore I shall state more clearly the reasons
why I feel that theorems like Euclid II 1--4 did not arise from
geometrical considerations.  \cite[pp.\;203--204]{Va76} (emphasis in
the original)
\end{quote}
Further details on van der Waerden's approach can be found in
Section~\ref{s25}.

We refrain from taking a position in the debate on
the narrow issue of geometric algebra as applied to Greek mathematics,
but point out that the debate has stimulated the articulation of
various approaches to the history of mathematics.  We will analyze how
the various approaches play themselves out in the study of the
pioneers of mathematical analysis including Fermat, Leibniz, Euler,
and Cauchy; see Section~\ref{s4}.  On whether other historians endorse
Unguru polarity, see Section~\ref{s7b}.

\subsection{Returning from escapades}
\label{s12}

Readers familiar with the tenor of the Unguru--Weil controversy will
not have been surprised by the tone of the ``break-in" remark found in
the 2001 book by Fried and Unguru (henceforth FU):
\begin{quote}
The mathematical and the historical approaches are antagonistic.
Whoever breaks and enters typically returns from his \emph{escapades}
with other spoils than the peaceful and courteous caller.
(Fried--Unguru \cite{Fr01}, 2001, p.\;406; emphasis on
\emph{escapades} added)
\end{quote}
For readers less familiar with the controversy, it may be prudent to
clarify that the unpeaceful and uncourteous caller allegedly involved
in the break-in is, to be sure, the mathematician, not the historian.

The break-in remark is duly reproduced in a recent essay in the
\emph{Journal of Humanistic Mathematics} by Unguru's disciple Fried
(\cite{Fr18}, 2018, p.\;7).  And yes, the ``escapades" are in the
original, both in FU and in Fried solo.

That Fried's mentor Unguru does not mince words with regard to Weil is
not difficult to ascertain.  Thus, one finds the following phrasing:
``Betrayals, Indignities, and Steamroller Historiography: Andr\'e Weil
and Euclid'' (Unguru \cite{Un18}, 2018, p.\;26, end of Section\;I).

Clearly, neither unpeaceful break-ins nor uncourteous escapades
represent \emph{legitimate} relationships to the mathematics of the
past.  The break-in remark makes the reader wonder about the precise
meaning of Fried's assurance that his plan is to catalog some of the
attitudes toward the history of mathematics ``\emph{without judgment}
as to whether they are necessarily correct or legitimate" (Fried
\cite[abstract, p.\;3]{Fr18}; emphasis added).

\subsection{Who is open-minded?}
\label{s13}

A further telling comment appears on the back cover of the FU book on
Apollonius of Perga:
\begin{quote}
Although this volume is intended primarily for historians of ancient
mathematics, its approach is fresh and engaging enough to be of
interest also to historians, philosophers, linguists, and
\emph{open-minded mathemati\-cians}.  (Fried--Unguru \cite{Fr01},
2001; emphasis added)
\end{quote}
Readers will not fail to notice that, of the four classes of scholars
mentioned, the mathematician is the only class limited by the
qualifier \emph{open-minded}; here FU don't appear to imply that a
mathematician is typically characterized by the qualifier.  Such a
polarizing approach (see further in Section~\ref{s52}) on the part of
FU is hardly consistent with Fried's professed idea of cataloguing
attitudes toward history ``without judgment as to whether they are
{\ldots} correct or legitimate.''  We will analyze the ideological
underpinnings of the FU approach in Section~\ref{s2}.

\section
{FU axiomatics: Discontinuity, \emph{tabula rasa}, antiplatonism}
\label{s2}

In this section, we will identify several axioms of historiography
according to the Unguru school.  In the Friedian scheme of things, it
is axiomatic that the proper view of the relation between the
mathematics of the past and that of the present is that of a
\emph{dis}continuity.  It is indeed possible to argue that it is (see
Section~\ref{s33} for an example in Fermat).  However, Fried appears
to take it for granted that a contrary view of \emph{continuity}
between past and present necessarily amounts to whig history or, more
politely, engaging in what Oakeshott described as a ``practical past"
and Grattan-Guinness \cite{Gr04} described as ``heritage" (see Fried
\cite{Fr18}, 2018, p.\;7).  Fried's attitude here is at odds with the
idea of historiography as seeking to ``explain why a certain
mathematical development happened'' (Fraser--Schroter, \cite{Fr19a},
2019, p.\;16).

\subsection{Axiom 1: Discontinuity}
\label{s31}

What emerges is the following axiom of Friedian historiographic
ideology:
\begin{quote}
\textbf{Axiom 1 (discontinuity).}  The proper attitude of a historian
toward the mathematics of the past is that of a discontinuity with the
mathematics of the present.
\end{quote}
We note that, while the discontinuity view may be appropriate in
certain cases, it is an assumption that needs to be argued rather than
posited as an axiom as in Fried.

Fried presents a taxonomy of various attitudes toward the history of
mathematics.  He makes it clear that he means to apply his taxonomy
rather broadly, and not merely to the historical work on Greek
mathematics:
\begin{quote}
[I]n most of the examples the mathematical past being considered by
one person or another is that of Greece. \ldots{} Nevertheless, as I
hope will be clear, the relationships evoked in the context of these
examples have little to do with the particular character of Greek
mathematics. \cite[p.\;8]{Fr18}
\end{quote}
We will examine the effectiveness of the FU approach in such a broader
context.

\subsection{Apollonius and \emph{mathematical conquerors}}
\label{s22}

Fried notes that his piece ``has been written in a light and playful
spirit" (ibid., p.\;8).  Such a light(-headed) spirit is reflected in
Fried's attitude toward Fermat.  Fermat's reconstruction of
Apollonius' \emph{Plane Loci} prompts Fried to place Fermat in a
category labeled ``mathematical conqueror" rather than that of a
historian.  The label also covers Descartes and Vi\`ete, whose ``sense
of the past has the unambiguous character of a `practical past', again
to use Oakeshott's term'' \cite[pp.\;11--12]{Fr18}.  Fried does not
spare the great magistrate the tedium of Mahoney's breezy journalese:
\begin{quote}
``Fermat was no antiquarian interested in a faithful reproduction of
Apollonius' original work; \ldots\;The Plane Loci was to serve as a
means to an end rather than an end in itself.''  (Mahoney as quoted by
Fried in \cite[p.\;11]{Fr18})
\end{quote}
Fried appears to endorse Mahoney's dismissive attitude toward Fermat's
historical work.  Regrettably, Fried ignores Recasens' more balanced
evaluation of Fermat's work on Apollonius' \emph{Plane Loci},
emphasizing its classical geometric style, and contrasting it with van
Schooten's:
\begin{quote}
Fermat's demonstration of Locus II-5 is presented in the classical
geometrical style of the day, though his conception was already
algebraic; [van] Schooten's is a pure exercise of analytic geometry.
(Recasens \cite{Re94}, 1994, p.\;315)
\end{quote}
What is the source of Fried's facile dismissal of Fermat's historical
scholarship?  While it is difficult to be certain, a clue is found in
the attitude of his mentor Unguru who wrote:
\begin{quote}
[Fermat] took the Greek problems away from their indigenous territory
into new and foreign lands.  Interestingly (and again, I think, quite
typically), Fermat did not see in his novel and revolutionary methods
strategies intrinsically alien to Greek mathematics, thus contributing
to the creation of the pervading and pernicious myth that there are
not indeed any substantive differences between the geometrical works
of the Greeks and the algebraic treatment of Greek mathematics by
post-Vi\`etan mathematicians.  (Unguru \cite{Un76}, 1976, p.\;775)
\end{quote}
Unguru describes Fermat's reading of Apollonius as contributing to a
``pernicious myth,'' and Fermat's reliance on Vieta's theory of
equations as ``alien'' to Greek mathematics; Fried apparently follows
suit.

By page 12, Fried takes on a group he labels
``mathematician-histo\-rians" whose fault is their interest in
historical continuity:
\begin{quote}
The mathematics of the past is still understood by them as continuous
with present mathematics. (Fried \cite[p.\;12]{Fr18})
\end{quote}
Again, adherence to the continuity view is cast without argument as a
fault (see Section~\ref{s31}).  Yet positing discontinuity as a
working hypothesis can make scholars myopic to important aspects of
the historical development of mathematics, and have a chilling effect
on attempts to explain why certain mathematical developments happened
(the Fraser--Schroter definition of the task of a historian); see
Section~\ref{s4} for some examples.

\subsection{Axiom 2: Tabula rasa}
\label{s3}

Unguru's opposite number Weil makes a predictable appearance in (Fried
\cite{Fr18}, 2018) on page 14, under the label \emph{privileged
observer}.  The label also covers Zeuthen and van der Waerden.  The
fault for this particular label is the desire to take advantage of
``their mathematical ideas \ldots\;to piece together the past"
(ibid.).  Fried's posture on Weil brings us to the next axiom of
Friedian historiographic ideology:
\begin{quote}
\textbf{Axiom 2 (tabula rasa).}  It is both possible and proper for
historians to refrain from using modern mathematical ideas.
\end{quote}
By page 16 we learn what authentic historians of mathematics do: they
\begin{quote}
take as their working assumption, a kind of \emph{null-hypothesis},
that there is a discontinuity between mathematical thought of the past
and that of the present.  \cite[p.\;16]{Fr18} (emphasis added)
\end{quote}
This formulation of Fried's discontinuity axiom (see
Section~\ref{s31}) has the advantage of being explicitly cast as a
\emph{hypothesis}.  Yet nothing about Fried's tone here suggests any
intention of actually exploring the \emph{validity} of such a
hypothesis.  A related objection to Fried's posture was raised by
Bl{\aa}sj\"o and Hogendijk (\cite{Bl18}, 2018, p.\;775), who argue
that ancient treatises may contain meanings and intentions that go
beyond the surface text, based on a study of Ptolemy's
\emph{Almagest}.  We will analyze Fried's Axiom\;2 in
Section~\ref{s4}.

\subsection{Axiom 3: Uprooting the Platonist deviation}
\label{s24}

\begin{figure}
\includegraphics[height=6in]{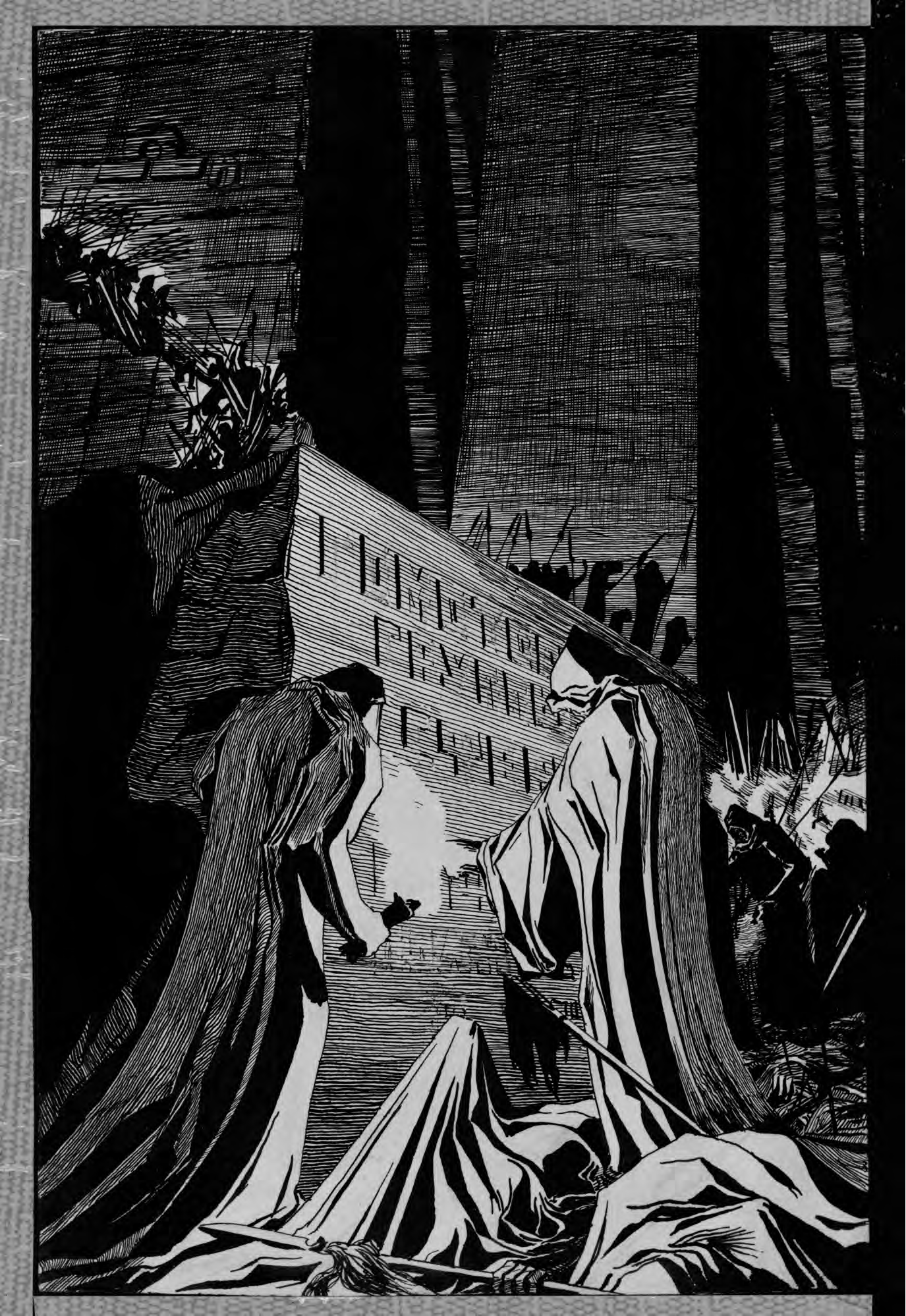}
\caption{\textsf{Mathematical Platonism.  A humorous illustration from
the book \emph{Homotopic Topology}.  The writing on the wall reads (in
Russian): ``Homotopy groups of spheres.''  Created by Professor
A. T. Fomenko, academician, Moscow State University.  Reproduced with
permission of the author.}}
\label{f1}
\end{figure}

The following additional axiom is discernible in the writings of
Unguru and his students.

\begin{quote}
\textbf{Axiom 3 (Mathematicians as Platonists)} Mathematicians
interested in history are predominantly Platonists; furthermore, their
beliefs (e.g., that mathematics is eternally true and unchanging)
interfere with their functioning as competent historians.
\end{quote}
Mathematicians interested in history are repeatedly described as
Platonists (see Figure~\ref{f1}) in the writings of Fried, Rowe, and
Unguru.  Thus, one finds the following comments (emphasis on
``Platonic'' added throughout):
\begin{enumerate}
[label={(Pl\,\theenumi)}]
\item
``[M]athematically minded historians {\ldots} assume tacitly or
explicitly that mathematical entities reside in the world of
\emph{Platonic} ideas where they wait patiently to be discovered by
the genius of the working mathematician" (Unguru--Rowe \cite{Un81},
1981, p.\;3, quoting \cite{Un79}).  Unguru and Rowe go on at length
(pages\;5 through 10) to attack van der Waerden's interpretation of
cunei\-form tablet BM 13901.%
\footnote{For a rebuttal of the Unguru--Rowe critique see Bl{\aa}sj\"o
(\cite{Bl16}, 2016, Section\;3.4).}
\item
``It has been argued that most contemporary historians of mathematics
are \emph{Platonists} in their approach.  They look in the past of
mathematics for the eternally true, the unchanging, the constant"
(Unguru--Rowe \cite{Un82}, 1982, p.\;47).  The problem with such an
approach is diagnosed as follows: ``If nothing changes there is no
history" (op.\;cit., p.\;48).
\item
``[T]he methodology embodied in `geometric algebra' {\ldots} is the
outgrowth of a \emph{Platonic} metaphysics that sees mathematical
ideas as disembodied beings, pure and untainted by any idiosyncratic
features'' (Fried--Unguru \cite{Fr01}, 2001, p.\;37).
\item
\label{pl4}
``That Apollonius was a skilled geometrical algebraist is clearly the
considerate opinion of Zeuthen.  It is an opinion based
\emph{exclusively} on a \emph{Platonic} philosophy of mathematics,
according to which one and the same mathematical idea remains the same
irrespective of its specific manifestations'' (Fried--Unguru
\cite{Fr01}, 2001, pp.\;47--48; emphasis added).
\item
``There is one mathematics, from its pre-historical beginnings to the
end of time, irrespective of its changing appearances over the
centuries.  This mathematics grows by accumulation and by a sharpening
of its standards of rigor, while, its rational, ideal, \emph{Platonic}
Kernel, remaining unaffected by the historical changes mathematics
undergoes, enjoys, as Hardy put it, immortality.  In short, proven
mathematical claims remain proven forever, no matter what the changes
are that mathematics is undergoing.  And since it is always possible
to present past mathematics in modern garb, ancient mathematical
accomplishments can be easily made to look modern and, therewith
seamlessly integrated into the growing body of mathematical
knowledge'' (Unguru \cite{Un18}, 2018, pp.\;19--20).
\item
``The \emph{Platonic} outlook embodied in Weil's statements, according
to which (1) mathematical entities reside in the world of
\emph{Platonic} ideas and (2) mathematical equivalence is tantamount
to historical equivalence, is inimical to history" (Unguru
\cite{Un18}, 2018, p.\;29).
\end{enumerate}
What is comical about this string of attempts to pin a
\emph{Platonist} label on scholars is the contrast between the extreme
care Unguru and his students advocate in working with primary
documents and sourcing every historical claim, on the one hand, and
the absence of such sources when it comes to criticizing scholars they
disagree with, on the other.  

This is not to say that mathematicians interested in history are never
Platonists.  Thus, in a recent volume by Dani--Papadopoulos, one
learns that
\begin{quote}
[T]hinkers in colonial towns in Asia Minor, Magna Graecia, and
mainland Greece, cultivated a love for systematizing phenomena on a
rational basis {\ldots} They appreciated purity, universality, a
certainty and an elegance of mathematics, the characteristics that all
other forms of knowledge do not possess.  (\cite{Pa19}, 2019, p.\;216)
\end{quote}

While such attitudes do exist, it remains that claiming your opponents
are Platonists without providing evidence is no more convincing than
claiming that Greek geometry had an algebraic foundation without
providing evidence, a fault Unguru and others impute to their
opponents.  Without engaging in wild-eyed accusations of Platonism
against scholars he disagreed with, Grattan-Guinness \cite{Gr96} was
able to enunciate a dignified objection to geometric algebra (for a
response see Bl{\aa}sj\"o \cite[Section\;3.10]{Bl16}).

With regard to Axiom\;3, it is worth noting that the broader the
spectrum of the culprits named by Unguru, the less credible his charge
of Platonism becomes.  Consider, for example, the claim in \ref{pl4}
above that Zeuthen's opinion is ``based \emph{exclusively} on a
Platonic philosophy of mathematics'' (emphasis added).  Unguru attacks
Heiberg in \cite[p.\;107]{Un75} with similar vehemence.  But how
credible would be a claim that the historiographic philosophy of the
\emph{philologist} Johan Ludvig Heiberg (of the \emph{Archimedes
Palimpsest} fame) is due to mathematical Platonist beliefs, especially
if no evidence is provided?

\subsection{Corry's universals}
\label{s25}

Axiomatizing tendencies similar to those of Fried, Rowe, and Unguru
manifest themselves in the writing by Unguru's student Corry, as well.

Engagement with Platonism and its discontents appears to be a constant
preoccupation in Corry's work.  He alludes to Platonism by using terms
as varied as ``eternal truth,'' ``essence of algebra,'' and
``universal properties.''  Thus, in 1997 he writes:
\begin{quote}
[Bourbaki] were extending in an unprecedented way the domain of
validity of the belief in the eternal character of mathematical
truths, from the body to the images of mathematical knowledge.  (Corry
\cite{Co97}, 1997, p.\;253)
\end{quote} 
In 2004 (originally published in 1996) he writes:
\begin{quote}
{\ldots} a common difficulty that has been manifest {\ldots} is the
attempt to define, by either of the sides involved, the ``essence" of
algebraic thinking throughout history.  Such an attempt appears, from
the perspective offered by the views advanced throughout the present
book, as misconceived.  (Corry \cite{Co04}, 2004, p.\;396).
\end{quote}
In 2013 we find:
\begin{quote}
{\ldots} the question about the ``essence of algebra'' as an
ahistorical category seems to me an ill-posed and uninteresting one.
(Corry \cite{Co13}, 2013, p.\;639)
\end{quote}
In 2007, Corry imputes to mathematicians a quest for ``universal
properties'' at the expense of historical authenticity.  He appears to
endorse Unguru's view of mathematicians as Platonists when he writes:
\begin{quote}
[I]n analyzing mathematics of the past mathematicians often look for
underlying mathematical concepts, regularities or affinities in order
to conclude about historical connection.  Mathematical affinity
necessarily follows from \emph{universal properties} of the entities
involved and this has often been taken to suggest a certain historical
scenario that `might be'.  But, Unguru warns us, one should be very
careful not to allow such mathematical arguments led [sic] us to
mistake historical truth (i.e., the `thing that \emph{has} been') with
what is no more than mathematically possible scenarios (i.e., the
`thing that \emph{might} be').  The former can only be found by direct
\emph{historical evidence}.  (Corry \cite{Co07}, 2007; emphasis on
``has'' and ``might'' in the original; emphasis on ``universal
properties'' and ``historical evidence'' added)
\end{quote}
Granted one needs direct historical evidence, as per Corry.  However,
where is the evidence that instead of looking for evidence,
mathematicians interested in history look for \emph{universal
properties}?  Such a view of mathematicians who are historians is
postulated axiomatically by Corry, similarly to Fried, Rowe, and
Unguru.  Corry goes on to claim that
\begin{quote}
[Unguru's 1975] work immediately attracted \emph{furious} reactions,
above all from three prominent mathematicians interested in the
history of mathematics: Andr\'e Weil, Bartel L. van der Waerden, and
Hans Freudenthal.  (ibid.; emphasis added)
\end{quote}
Corry's remark is specifically characterized by the attitude of
looking for ``the thing that might be" rather than the ``thing that
has been", a distinction he mentions in the passage quoted above.
Namely, once Corry postulates a \emph{universal}-seeking attitude on
the part of Weil, van der Waerden, and Freudental, it then naturally
follows, for Corry, that they would necessarily react ``furiously" to
Unguru's work.  Corry does present evidence of presentist attitudes in
historiography in connection with interpreting the Pythagorean
discovery of the incommensurability of the diagonal and the side of
the square.  However, Corry presents evidence of such attitudes not in
the writings of Weil, van der Waerden, or Freudental, but rather those
of{\ldots} Carl Boyer (\cite{Bo68}, 1968, p.\;80).%
\footnote{In \cite{Co07}, Corry criticizes attempts to deduce a purely
historical claim merely from ``underlying mathematical affinity.''
Corry provides the following example: ``It is thus inferred that the
Pythagoreans proved the incommensurability of the diagonal of a square
with its side exactly as we nowadays prove that $\sqrt{2}$ is an
irrational number.'' Corry's example is followed by a reference to
(Boyer \cite{Bo68}, 1968, p.\;80).  On that page, Boyer wrote: ``A
third explanation [of the expulsion of Hippasus from the Pythagorean
brotherhood] holds that the expulsion was coupled with the disclosure
of a mathematical discovery of devastating significance for
Pythagorean philosophy---the existence of incommensurable
magnitudes.''}
For a discussion of the shortcomings of Boyer's historiographic
approach see Section~\ref{s43}.  Note that Weil is just as sceptical
as Corry about claims being made on behalf of the Pythagoreans; see
(Weil \cite{We84}, 1984, pp.\;5, 8).  We will analyze Corry's
problematic criticisms of van der Waerden, Freudenthal, and Weil
respectively in Sections~\ref{s26}, \ref{s27}, and \ref{s28}.

\subsection
{van der Waerden on Diophantus and Arabic algebra}
\label{s26}

Contrary to the claims emanating from the Unguru school, some of
Unguru's opponents specifically denied being Platonist.  Thus, van der
Waerden wrote:
\begin{quote}
I am simply not a Platonist.  For me mathematics is not a
contemplation of essences but intellectual construction.  (van der
Waerden as translated in Schappacher \cite{Sc07}, 2007, p.\;245)
\end{quote}
It is instructive to contrast Unguru's attitude toward van der Waerden
with Szab\'o's.  Szab\'o's book (\cite{Sz78}, 1978) deals with van der
Waerden at length, but there is no trace there of any allegation of
Platonist deviation.  On the contrary, Szab\'o's book relies on van
der Waerden's historical scholarship, as noted also by Folkerts
(\cite{Fo69}, 1969).  The book does criticize van der Waerden for what
Szab\'o claims to be an over-reliance on translations.  Szab\'o
discusses this issue in detail in the context of an analysis of the
meaning of the Greek
term~$\delta\acute\upsilon\nu\alpha\mu\iota\varsigma$
(\emph{dynamis}),%
\footnote{At the risk of committing precisely the type of inaccuracy
criticized by Szab\'o, one could translate \emph{dynamis} roughly as
``the squaring operation''.}
as mentioned by Folkerts.

In reality, van der Waerden's 1976 article contains no sign of the
``fury'' claimed by Corry (see Section~\ref{s25}).  Perhaps the most
agitated passage there is van der Waerden's rebuttal of a spurious
claim by Unguru (which is echoed thirty years later by Corry in
\cite{Co07}):
\begin{quote}
We (Zeuthen and his followers) feel that the Greeks started with
algebraic problems and translated them into geometric language.
Unguru thinks that we argued like this: We found that the theorems of
Euclid II can be translated into modern algebraic formalism, and that
they are easier to understand if thus translated, and this we took as
`the proof that this is what the ancient mathematician had in mind'.
Of course, this is nonsense.  We are not so weak in logical thinking!
The fact that a theorem can be translated into another notation does
not prove a thing about what the author of the theorem had in mind.
(van der Waerden \cite{Va76}, 1976, p.\;203)
\end{quote}
What one \emph{does} find in van der Waerden's article is a specific
rebuttal of Corry's claim concerning an alleged search for
\emph{universal properties} (see Section~\ref{s25}).  Here van der
Waerden was responding to Unguru, who had claimed that algebraic
thinking involves
\begin{quote}
Freedom from any ontological questions and commitments and, connected
with this, abstractness rather than intuitiveness.  (Unguru
\cite{Un75}, 1975, p.\;77)
\end{quote}
Van der Waerden responded by rejecting Unguru's characterisation of
the algebra involved in his work on Greek mathematics, and pointed out
that what he is referring to is
\begin{quote}
algebra in the sense of Al-Khw\={a}rizm\={\i}, or in the sense of
Cardano's `Ars magna', or in the sense of our school algebra.  (van
der Waerden \cite{Va76}, 1976, p.\;199)
\end{quote}
Thus, van der Waerden specifically endorsed a similarity between
Arabic premodern algebra and Greek mathematics, and analyzed
Diophantus specifically in \cite[p.\;210]{Va76}.

A similarity between Arabic premodern algebra and the work of
Diophantus is also emphasized by Christianidis (\cite{Ch18}, 2018).
The continuation of the passage from van der Waerden is more
problematic from the point of view of \cite{Ch18}.  Here he wrote:
\begin{quote}
Algebra, then, is: the art of handling algebraic expressions
like~$(a+b)^2$ and of solving equations like \mbox{$x^2+ax=b$}.  (van
der Waerden \cite{Va76}, 1976, p.\;199)
\end{quote}
The viewpoint expressed here is at odds with the emphasis in
\cite{Ch18} on the fact that premodern algebra did not deal with
equations, polynomials were not sums but rather aggregates, and the
operations stipulated in the problem were performed before the
statement of the equation.  However, apart from these important
points, van der Waerden's notion of Greek mathematics as close to
Arabic premodern algebra is kindred to the viewpoint elaborated in
\cite{Ch18} and \cite{Ch19}.%
\footnote{It would be more difficult to bridge the gap between the
positions of Weil and \cite{Ch18}, since Weil claims that ``there is
much, in Diophantus and in Vi\`ete's \emph{Zetetica}, which in our
view pertains to algebraic geometry'' \cite[p.\;25]{We84}, whereas
\cite{Ch18} specifically distances itself from attempts to interpret
Diophantus in terms of algebraic geometry.}

The article \cite{Ch18} elaborates a distinction between \emph{modern
algebra} and \emph{premodern algebra}.  The latter term covers both
Arabic sources and Diophantus.  The position presented in \cite{Ch18}
is clearly at odds with Unguru, who wrote:
\begin{quote}
With Vi\`ete algebra becomes the very language of mathematics; in
Diophantus' \emph{Arithmetica}, on the other hand, we possess merely a
refined auxiliary tool for the solution of arithmetical problems
{\ldots} (Unguru \cite{Un75}, 1975, p.\;111, note\;138)
\end{quote}

Regardless of how close the positions of van der Waerden and
\cite{Ch18} can be considered to be, there is no mention of
\emph{universal properties} in van der Waerden's work on Greek
mathematics.  The \emph{universals} appear to be all Corry's, not van
der Waerden's.

\subsection{Corry's shift on Freudenthal}
\label{s27}

Similarly instructive is Corry's--as we argue--variable position on
Freu\-denthal in connection with Platonism and Bourbaki.  The standard
story on Bourbaki is the one of mathematical Formalism and
\emph{structures}.  In Corry's view, there is some question concerning
how different this is, in Bourbaki's case, from Platonism.  Here is
what Corry wrote in his book in 1996:
\begin{quote} 
The above-described mixture [in Bourbaki] of a declared formalist
philosophy with a heavy dose of \emph{actual Platonic belief} is
illuminating in this regard.  The formalist imperative, derived from
that ambiguous position, provides the necessary background against
which Bourbaki's drive to define the formal concept of
\emph{structure} and to develop some immediate results connected with
it can be conceived.  The \emph{Platonic} stand, on the other hand,
which reflects Bourbaki's true working habits and beliefs, has led the
very members of the group to consider this kind of conventional,
formal effort as superfluous.  (Corry \cite{Co96}, 1996, p.\;311;
emphasis on ``structure'' in the original; emphasis on ``actual
Platonic belief'' and ``Platonic'' added)
\end{quote}
On page 336 in the same book, Corry quotes Freudenthal's biting
criticism of Piaget's reliance on Bourbaki and their concept of
structure as an organizing principle:
\begin{quote}
The most spectacular example of organizing mathematics is, of course,
Bourbaki.  How convincing this organization of mathematics is!  So
convincing that Piaget could rediscover Bourbaki's system in
developmental psychology.  {\ldots} Piaget is not a mathematician, so
he could not know how unreliable mathematical system builders are.
(Freudenthal as quoted in Corry \cite{Co96}, 1996, p.\;336).
\end{quote}
The same passage is quoted in the earlier article (Corry \cite{Co92},
1992, p.\;341).

The index in Corry's book on Bourbaki contains an ample supply of
entries containing the term \emph{universal}, including
\emph{universal constructs} and \emph{universal problems}, a constant
preoccupation of Bourbaki's which can also be seen as a function of
their Platonist background philosophy in the sense Corry outlined in
\cite[p.\;311]{Co96}, where Corry speaks of Bourbaki's ``actual
Platonic belief.''

There is a clear contrast, in Corry's mind, between Platonism,
universals, Bourbaki, and Piaget, on one side of the debate, and
Freudenthal with his clear opposition to both Bourbaki and Piaget, on
the other.  Freudenthal's opposition tends to undercut the idea of
Freudenthal as Platonist, which in any case is at odds with
Freudenthal's pragmatic position on mathematics education; see e.g.,
La Bastide (\cite{La15}, 2015).

On the other hand, Corry's article (\cite{Co07}, 2007) includes
Freudenthal on the list of the Platonist mathematical culprits (van
der Waerden, Freudenthal, Weil) that has been made standard by Unguru.
According to Corry 2007, these scholars are in search of
\emph{universal properties}; see Section\;\ref{s25}.  This fits with
Unguru's take on Freudenthal, but is at odds with what Corry himself
wrote about Freudenthal a decade earlier, as documented above.

In fact, Freudenthal specifically sought to distance himself from
Platonism in (\cite{Fr78}, 1978, p.\;7).  Freudenthal's interest in
Intuitionism is discussed in \cite[p.\;42]{La15}.  He published at
least two papers in the area: \cite{Fr37a} and \cite{Fr37b}.  This
interest similarly points away from Platonism, contrary to Corry's
claim.

%

\subsection{Weil: internalist or externalist?}
\label{s28}

Corry attacks both Weil and Bourbaki as Platonist, and dismisses
Bourbaki's volume on the history of mathematics \cite{Bo94} as
``royal-road-to-me" historiography, in (Corry \cite{Co07}, 2007).
Paumier and Aubin make a more specific claim against the Bourbaki
volume generally and Weil's historiography in particular.  Namely,
they refer to the volume as ``internalist history of concepts''
\cite[p.\;185]{Pa16}, and imply that the same criticism applies to
Weil's historiography, as well.  In a related vein, Kutrov\'atz casts
Unguru and Szab\'o as externalists and Weil as internalist and
Platonist in \cite{Ku02}.

To evaluate such criticisms of Weil, the Bourbaki volume is of limited
utility since it was of joint authorship.  We will examine instead
Weil's own book \emph{Number theory.  An approach through history.
From Hammurapi to Legendre} (\cite{We84}, 1984).  Does the
\emph{internalist} criticism apply here?  

To answer the question, we would need to agree first on the meaning of
\emph{internalist} and \emph{externalist}.  If we posit that
historical work is \emph{externalist} if it is written by Unguru, his
disciples, and their cronies, then there is little hope for Weil.
There is perhaps hope with a less partisan definition, such as
``historiography that takes into account the contingent details of the
historical period and its social context, etc.''
%
%
It is clear that historical and social factors are important.  For
instance, one obtains a distorted picture of the mathematics of
Gregory, Fermat, and Leibniz if one disregards the fierce religious
debates of the 17th century (see references listed in
Section~\ref{s41b}).

Now it so happens that Weil's book \cite{We84} does contain detailed
discussions of the historical context.  Weil's book is not without its
shortcomings.  For instance, when Weil mentions that Bachet
``extracted from Diophantus the conjecture that every integer is a sum
of four squares, and asked for a proof'' \cite[p.\;34]{We84}, the
reader may well feel disappointed by the ambiguity of the verb
``extracted'' and the absence of references.  However, what interests
us here is the validity or otherwise of the contention (implicit in
Unguru and Corry and explicit in Paumier--Aubin and Kutrov\'atz) that
Weil was \emph{internalist}.  Was Weil internalist as charged?

Weil mentions, for instance, that Euler was first motivated to look at
the problem of Fermat primes $2^{2^n}+1$ by his correspondent Goldbach
\cite[p.\;172]{We84}.  To give another example, Weil mentions that
Fermat learned Vieta's symbolic algebra through his visits to
d'Espagnet's private library in Bordeaux in the 1620s
\cite[p.\;39]{We84}.  Such visits took place many years before Vieta's
works were published in 1646 by van Schooten.  In particular, the
Fermat--d'Espagnet contact was instrumental in Fermat's formulation of
his method of adequality (see Section~\ref{s32}) relying as it did on
Vieta's symbolic algebra.  Such examples undermine the Paumier--Aubin
claims, such as the following:
\begin{enumerate}
\item
the charge of ``an `internalist history of concepts' which has only
little to say about the way in which mathematics emerged from the
interaction of groups of people in specific circumstances''
\cite[p.\;187]{Pa16};
\item
the claim that ``The focus on ideas erased much of the social dynamics
at play in the historical development of mathematics''
\cite[p.\;204]{Pa16}.
\end{enumerate}
As we showed, Weil does take the interactions and the dynamics into
account.

\section{Some case studies}
\label{s4}

We identified Fried's \emph{tabula rasa} axiom in Section~\ref{s3},
and will analyze it in more detail in this section.  It seems that
while the axiom may be appropriate in certain cases, it is an
assumption that needs to be argued rather than merely postulated.
Such a need to argue the case applies to the very possibility itself
of a ``tabula rasa" attitude in the first place:
\begin{enumerate}
 [label={(TR\theenumi)}]
\item
\label{i1}
Can historians of mathematics truly view the past without the lens of
modern mathematics?
\item
\label{i2}
Have historians been successful in such an endeavor?
\end{enumerate}

Whereas it may be difficult to rule out the theoretical possibility of
an affirmative answer to \ref{i1}, a number of recent studies suggest
that in practice, the answer to \ref{i2} is often negative, as we will
discuss in Section~\ref{s41b}.

\subsection{History of analysis}
\label{s41b}

Some historians of 17th through 19th century mathematical analysis,
while claiming to reject insights provided by modern mathematics in
their interpretations, turn out themselves to be \emph{privileged
observers} in Fried's sense (see Section~\ref{s3}) though still in
denial.  Namely, they operate within a conceptual scheme dominated by
the mathematical framework developed by Weierstrass at the end of the
19th century, as argued in recent studies in the following cases:
\begin{itemize}
\item
Fermat, in Katz et al.\;(\cite{13e}, 2013) and Bair et
al.\;(\cite{18d}, 2018);
\item
Gregory, in Bascelli et al.\;(\cite{18f}, 2018);
\item
Leibniz, in Sherry--Katz (\cite{14c}, 2014), Bascelli et
al.\;(\cite{16a}, 2016), Bl{\aa}sj\"o (\cite{Bl17}, 2017), and Bair et
al.\;(\cite{18a}, 2018);
\item
Euler, in Kanovei et al.\;(\cite{15b}, 2015) and Bair et
al.\;(\cite{17b}, 2017);
\item
Cauchy, in Bair et al.\;(\cite{17a}, 2017); Bascelli et
 al.\;(\cite{18e}, 2018); and Bair et al.\;(\cite{19a}, 2019 and
 \cite{20a}, 2020).
\end{itemize}
The pattern that emerges from these studies is that some modern
historians, limited in their knowledge of modern mathematics, tend to
take a narrow view, that in some cases borders on naivete, of the work
of the great mathematical pioneers of the 17--19th centuries (see
Sections~\ref{s44} and \ref{s45} for examples).

The Fermat historian Mahoney is a case in point.  Weil pointed out
numerous historical, philological, and mathematical errors in
Mahoney's work on Fermat; see \cite{We73}.  Yet in the Friedian scheme
of things, Weil is neatly shelved away on the \emph{privileged
observer} shelf, whereas Mahoney's work, breezy journalese and all, is
blithely assumed to reside in that rarefied stratum called authentic
history of mathematics, and relied upon to pass judgments on the value
of the historical work by the great Pierre de Fermat (see
Section~\ref{s31}).

\subsection{Fermat's adequality}
\label{s32}

Fermat used the method of \emph{adequality} to find maxima and minima,
tangents, and solve other problems.

To illustrate Fermat's method, consider the first example appearing in
his \emph{Oeuvres} \cite[p.\;134]{Ta1}.  Fermat considers a segment of
length~$B$, splits it into variable segments of length~$A$ and~$B-A$,
and seeks to maximize the product~$A(B-A)$, i.e.,
\begin{equation}
\label{e31}
BA-A^2.
\end{equation}
Next, Fermat replaces~$A$ by~$A+E$ (and~$B-A$ by~$B-A-E$).  There is a
controversy in the literature as to exact nature of Fermat's~$E$, but
for the purposes of following the mathematics it may be helpful to
think of~$E$ as small.  Fermat goes on to expand the corresponding
product as follows:
\begin{equation}
\label{e32}
BA-A^2+BE-2AE-E^2.
\end{equation}
In order to compare the expressions \eqref{e31} and \eqref{e32},
Fermat removes the terms independent of~$E$ from both expressions, and
forms the relation
\begin{equation}
\label{e33}
BE\adequal 2AE+E^2,
\end{equation}
also referred to as \emph{adequality}.  In the original, the term
\emph{ad{\ae}quabitur} appears where we used the symbol~$\adequal$.%
\footnote{A symbol similar to $\adequal$ was used several decades
later by Leibniz, interchangeably with~$=$, to denote a relation of
generalized equality.}
We will present the final part of Fermat's solution in
Section~\ref{s34}.

\subsection{van Maanen's summary}

Fermat's method is described as follows by van Maanen:%
\footnote{In place of Fermat's~$E$, van Maanen uses a lower-case~$e$.
The pieces of notation~$I(x)$,~$x_M$, and~$=$ are van Maanen's.}
\begin{quote}
Fermat seems to have based his method for finding a maximum or minimum
for a certain algebraic expression~$I(x)$ on a double root argument,
but in practice the algorithm was used in the following slightly
different form.  Fermat argued that if the extreme value is attained
at~$x_M$,~$I(x)$ is constant in an infinitely small neighborhood%
\footnote{Describing Fermat's method in terms of the infinitely small
is not entirely uncontroversial and is subject to debate; for details
see Bair et al.\;(\cite{18d}, 2018).}
of~$x_M$.  Thus, if~$E$ is very small,~$x_M$ satisfies the equation
$I(x+E)=I(x)$.  \cite[p.\;52]{Va03}
\end{quote}
Fermat never actually formed an algebraic relation (using Vieta's
symbolic algebra) of \emph{ad{\ae}quabitur} between the
expressions~$I(x+E)$ and~$I(x)$.  The kind of relation he did form is
illustrated in formula~\eqref{e33} in Section~\ref{s32}.  Van Maanen
provides the following additional explanations:
\begin{quote}
This expression states that close to the extreme value,%
\footnote{We added the comma for clarity.}
lines parallel to the~$x$-axis will intersect the graph of~$I$ in two
different points, but the extreme is characterised by the fact that
these {\ldots} parallels turn into the tangent [line] and the points
of intersections reduce to one point which counts twice.  The common
terms in~$x$ are removed from the equation~$I(x+E)=I(x)$ and the
resulting equation divided by~$E$.  Any remaining terms are deleted,
and~$x_M$ is solved from the resulting equation.  (ibid.)
\end{quote}
While the summary by van Maanen does not mention the possibility of
dividing by $E^2$, it is important to note that in Fermat's
descriptions of the method, Fermat does envision the possibility of
dividing by higher powers of~$E$ in the process of obtaining the
extremum.

\subsection{Squaring both sides}
\label{s34}

In the example presented in Section~\ref{s32}, the term~$BE$ and the
sum~$2AE+E^2$ originally both appeared in the expression~\eqref{e32},
but appear on different sides in relation~\eqref{e33} (all with
positive sign).  The remainder of Fermat's algorithm is more familiar
to the modern reader: one divides both sides by~$E$ to produce the
relation~$B\adequal 2A+E$, and discards the summand~$E$ to obtain the
solution~$A=\frac{B}{2}$.

For future reference, we note that a relation of type~\eqref{e33} can
be squared to produce a relation of type~$(BE)^2 \adequal
(2AE+E^2)^2$.  In this particular example, the relation need not be
squared.  However, in an example involving square roots one needs to
square both sides at a certain stage to eliminate the radicals; see
Section~\ref{s33}.  Meanwhile, once one passes to the difference
$I(x+E)-I(x)$ (to use van Maanen's notation), such an opportunity is
lost.

Fermat never performed the step of carrying all the terms to the
left-hand side of the relation so as to form the difference
$I(x+E)-I(x)$;%
\footnote{Fermat historian Breger did in (\cite{Br13}, 2013, p.\;27);
for details see Bair et al.\;(\cite{18d}, 2018, Section\;2.6,
p.\;573).}
nor did Fermat ever form the quotient~$\frac{I(x+E)-I(x)}{E}$ familiar
to the modern reader.  In Section~\ref{s33} we will compare the
treatment of this aspect of Fermat's method by a historian and a
mathematician.

\subsection{Experiencing~$E^2$}
\label{s33}

The perspective of Unguruan polarity can lead historians to devote
insufficient attention to the actual mathematical details and
ultimately to historical error.  Thus, Mahoney claimed the following:
\begin{quote}
In fact, in the problems Fermat worked out, the proviso of repeated
division by~$y$ [i.e.,~$E$] was unnecessary.  But, thinking in terms
of the theory of equations, Fermat could imagine, \emph{even if he had
not experienced}, cases in which the adequated expressions contained
nothing less than higher powers of~$y$.  (Mahoney \cite{Ma94}, 1994,
p.\;165; emphasis added)
\end{quote}
Mahoney assumed that Fermat ``had not experienced'' cases where
division by~$E^2$ was necessary.  Meanwhile, Giusti analyzes an
example ``experienced'' by Fermat which involves radicals, and which
indeed leads to division by~$E^2$.  The example (Fermat \cite{Ta1},
p.\;153) involves finding the maximum of the expression
$A+\sqrt{BA-A^2}$ (here~$B$ is fixed).  In the process of solution, a
suitable relation of adequality, as in formula~\eqref{e33}, indeed
needs to be \emph{squared} (see Section~\ref{s34}).  Giusti concludes:
\begin{quote}
Ce qui nous int\'eresse dans ce cas est qu'il donne en exemple une
ad\'egalit\'e o\`u les termes d'ordre le plus bas sont en~$E^2$.
Comme on sait, dans l'\'enonciation de sa r\`egle Fermat parlait de
division par~$E$ ou par une puissance de~$E$ {\ldots} Plusieurs
commentateurs ont soutenu {\ldots} que Fermat avait commis ici une
\emph{erreur} {\ldots} On doit donc penser que dans un premier moment
Fermat avait trait\'e les quantit\'es contenant des racines avec la
m\'ethode usuelle, qui conduisait parfois \`aà la disparition des
termes en~$E$, et qui ait tenu compte de cette \'event\-ualit\'e dans
l'\'enonciation de la r\`egle g\'en\'erale.  (Giusti \cite{Gi09},
2009, Section\;6, pp.\;75--76; emphasis added).
\end{quote}
What Giusti is pointing out is that in this particular application of
adequality in a case involving radicals, division by~$E^2$ (and not
merely by~$E$) is required.  Thus, the error is Mahoney's, not
Fermat's.

A first-rate analyst and differential geometer, Giusti was able to
appreciate the \emph{discontinuity} between Fermat's method of
adequality, on the one hand, and the modern~$\frac{I(x+E)-I(x)}{E}$
perspective, on the other, better than many a Fermat historian.  More
generally, a scholar's work should be evaluated on the basis of its
own merits rather than which class he primarily belongs to, be it
historian, mathematician, or philologist.

Appreciating discontinuity is not the prerogative of Unguru's adepts,
contrary to strawman accounts found in Unguru (\cite{Un18}, 2018) and
Guicciardini (\cite{Gu18}, 2018).  The portrait of a mathematician's
view of his discipline dominated by mathematical Platonism as found in
Unguru and his students (as detailed in Section~\ref{s24}) as well as
Guicciardini is similarly a strawman caricature, as when Guicciardini
elaborates on ``the perfect embodiment of the immutable laws of
mathematics written in the sky for eternity" \cite[p.\;148]{Gu18} and
claims that ``[t]he mathematician's world is the world of Urania''
(op.\;cit., p.\;150).

\subsection
{Why certain developments happened: Euler to Cauchy}
\label{s43}

Analyzing the differences between 18th and 19th century analysis,
Fraser and Schroter observe:
\begin{quote}
The decline of [Euler's] formalism stemmed mainly from its limitations
as a means of generating useful results.  Moreover, as methods began
to change, an awareness of formalism's apparent difficulties and even
contradictions lent momentum to efforts to rein it in. 
(Fraser--Schroter \cite{Fr20}, 2020, Section 3.3)
\end{quote}
Note that Fraser and Schroter are analyzing Euler's own work itself
here, rather than its reception by Cauchy.  Fraser and Schroter
continue:
\begin{quote}
Euler had been confident that the ``out-there'' objectivity of algebra
secured the generality of his formal techniques, but Cauchy demanded
that generality be found within mathematical methods themselves. In
his [textbook] \emph{Cours d'analyse} of 1821 Cauchy rejected
formalism in favour of a fully quantitative analysis.  (ibid.)
\end{quote}
Fraser and Schroter feel that the limitations and difficulties of
Euler's variety of algebraic formalism can be fruitfully analyzed from
the standpoint of considerably later developments, notably Cauchy's
``quantitative analysis.''  In their view, it is possible to comment
on the shortcomings of Euler's algebraic formalism and the reasons for
this particular development from Euler to Cauchy without running the
risk of anachronism.  Meawhile, it is clear that the Fraser--Schroter
approach may run afoul of both the discontinuity axiom (see
Section~\ref{s31}) and the tabula rasa axiom (see Section~\ref{s3}).

The issue of anachronism was perceptively analyzed by Ian Hacking
(\cite{Ha14}, 2014) in terms of the distinction between the butterfly
model and the Latin model for the development of a scientific
discipline.  Hacking contrasts a model of a deterministic (genetically
determined) biological development of animals like butterflies (the
egg--larva--cocoon--butterfly sequence), with a model of a contingent
historical evolution of languages like Latin.  Emphasizing determinism
over contingency can easily lead to anachronism; for more details see
Bair et al.\;(\cite{19a}, 2019).

Similarly to Hacking, Fraser notes the danger for a historian in the
adoption of a model based on an analogy with the pre-determined
evolution of a biological organism.  In his review of Boyer's book
\emph{The concepts of the calculus}, Fraser comments on the risks of
anachronism:
\begin{quote}
[Boyer's] focus on the development of concepts through time may
reflect as well an embrace of the metaphor of a plant or animal
organism.  The concept undergoes a progressive development, moving in
a directed and pre-determined way from its origins to an adult and
completed form.  {\ldots} The possibility of introducing anachronisms
is almost inevitable in such an approach, and to a certain degree this
is true of Boyer's book.  (Fraser \cite{Fr19}, 2019, p.\;18)
\end{quote}
Fraser specifically singles out for criticism Boyer's teleological
view of mathematical analysis as inexorably progressing toward the
ultimate \emph{Epsilontik} achievement:
\begin{quote}
[Boyer] seemed to view the eighteenth-century work as exploratory or
approximative as the subject moved inexorably in the direction of the
arithmetical limit-based approach of Augustin-Louis Cauchy and Karl
Weierstrass.  (op. cit., p.\;19)
\end{quote}
We will report on two additional cases of such teleological thinking
in the historiography of mathematics in Sections~\ref{s44} and
\ref{s45}.

\subsection{Leibnizian infinitesimals}
\label{s44}

Boyer-style, \emph{Epsilontik}-oriented teleological readings of the
history of analysis (see Section~\ref{s43}) are common in the
literature.  Thus, Ishiguro interprets Leibnizian infinitesimals as
follows:
\begin{quote}
It seems that when we make reference to infinitesimals in a
proposition, we are not designating a fixed magnitude incomparably
smaller than our ordinary magnitudes. Leibniz is saying that whatever
small magnitude an opponent may present, one can assert the existence
of a smaller magnitude. In other words, we can paraphrase the
proposition with a universal proposition with an embedded existential
claim.  (Ishiguro \cite{Is90}, 1990, p.\;87)
\end{quote}
What is posited here is the contention that when Leibniz wrote that
his incomparable (or inassignable)~$dx$, or~$\epsilon$, was smaller
than every given (assignable) quantity~$Q$, what he really meant was
an alternating-quantifier clause (universal quantifier~$\forall$
followed by an existential one~$\exists$) to the effect that for each
given~$Q>0$ there exists an~$\epsilon>0$ such that~$\epsilon<Q$.  Such
a logical sleight of hand goes under the name of the
\emph{syncategorematic interpretation}.  Here the author is
interpreting Leibniz as thinking like Weierstrass (see also
Section~\ref{sync}).  For details see Bascelli et al.\;(\cite{16a},
2016) and Bair et al.\;(\cite{18a}, 2018).

\subsection{Cauchyan infinitesimals}
\label{s45}

In a similar vein, Siegmund-Schultze views Cauchy's use of
infinitesimals as a step \emph{backward}:
\begin{quote} 
There has been \ldots{}\;an intense historical discussion in the last
four decades or so how to interpret certain apparent remnants of the
past or -- as compared to J. L. Lagrange's (1736--1813) rigorous
`Algebraic Analysis' -- even \emph{steps backwards} in Cauchy's book,
particularly his use of infinitesimals {\ldots} (Siegmund-Schultze
\cite{Si09}, 2009; emphasis added)
\end{quote}
Siegmund-Schultze's reader will have little trouble reconstructing
exactly which direction a \emph{step forward} may have been in.
Grabiner similarly reads Cauchy as thinking like Weierstrass; for
details see Bair et al.\;(\cite{19a}, 2019).

\subsection{History, heritage, or escapade?}

Significantly, in his essay Fried fails to mention the seminal
scholarship of Reviel Netz on ancient Greek mathematics (see e.g.,
\cite{Ne02b}, \cite{Ne01}, \cite{Ne02}).  Would, for example, Netz's
detection of traces of infinitesimals in the work of Archimedes be
listed under the label of history, heritage, or ``escapade" (to quote
Fried)?  Would an argument to the effect that the \emph{procedures}
(see Section~\ref{sync}) of the Leibnizian calculus find better
proxies in modern infinitesimal frameworks than in late 19th century
Weierstrassian ones, rank as history, heritage, or escapade?  Would an
argument to the effect that Cauchy's definition of continuity via
infinitesimals finds better proxies in modern infinitesimal frameworks
than in late 19th century Weierstrassian ones, rank as history,
heritage, or escapade?  Unfortunately, there is little in Fried's
essay that would allow one to explore such questions.

\subsection{Procedures vs ontology}
\label{sync}

The procedures/ontology distinction elaborated in B{\l}aszczyk et al.
(\cite{17d}, 2018) can be thought of as a refinement of
Grattan-Guinness' history/heritage distinction.  Consider for instance
Leibnizian infinitesimal calculus.  Without the procedure/ontology
distinction, interpreting Leibnizian infinitesimals in terms of modern
infinitesimals will be predictably criticized for utilizing history as
\emph{heritage}.  What some historians do not appreciate sufficiently
is that, in an ontological sense, interpreting Leibniz in
Weierstrassian terms is just as much heritage.  Surely talking about
Leibniz in terms of ultrafilters%
\footnote{See e.g., Fletcher et al.\;\cite{17f} for a technical
explanation.}
is not writing history; however, analyzing Leibnizian
\emph{procedures} in terms of those of Robinson's procedures is better
history than a lot of what is written on Leibniz by received
historians and philosophers (who have pursued a syncategorematic
reading of Leibnizian infinitesimals; see Section~\ref{s44}), such as
Ishiguro, Arthur, Rabouin, and others; for details see Bair et
al.\;(\cite{18a}, 2018).

We summarize some of the arguments involved.  \textbf{1.}\;Leibniz
made it clear on more than one occasion that his infinitesimals
violate Euclid Definition\;V.5 (Euclid\;V.4 in modern editions), which
is a version of what is known today as the Archimedean axiom; see
e.g., (Leibniz \cite{Le95b}, 1695, p.\;322).  In this sense, the
procedures in Leibniz are closer to those in Robinson than those in
Weierstrass.  \textbf{2.}\;If one follows Unguru's strictures and
Fried's tabula rasa, one can't exploit \emph{any} modern framework to
interpret Leibniz; however, in practice the syncategorematic society
interpret Leibniz in Weierstrassian terms, so the Unguruan objection
is a moot point as far as the current debate over the Leibnizian
calculus is concerned.  \textbf{3.}\;While modern foundations of
mathematics were clearly not known to Leibniz, it is worth pointing
out that this applies both to the set-theoretic foundational
\emph{ontology} of the classical Archimedean track, and to Robinson's
non-Archimedean track.  But as far as Leibniz's \emph{procedures} are
concerned, they find closer proxies in Robinson's framework than in a
Weierstrassian one.  For example, Leibniz's law of continuity is more
readily understood in terms of Robinson's transfer principle than in
any Archimedean terms.  \textbf{4.}\;The syncategorematic society
seems to experience no inhibitions about interpreting Leibnizian
infinitesimals in terms of alternating quantifiers (see
Section~\ref{s44}), which are conspicuously absent in Leibniz himself.
Meanwhile, Robinson's framework enables one to interpret them without
alternating quantifiers in a way closer to Leibniz's own procedures.
\textbf{5.}\;On several occasions Leibniz mentions a distinction
between \emph{inassignable} numbers like~$dx$ or~$dy$, and (ordinary)
\emph{assignable} numbers; see e.g., his \emph{Cum Prodiisset}
\cite{Le01c} and \emph{Puisque des personnes}\ldots\;\cite{Le05b}.
The distinction has no analog in a traditional Weierstrassian
framework.  Meanwhile, there is a ready analog of standard and
nonstandard numbers, either in Robinson's \cite{Ro66} or in Nelson's
\cite{Ne77} framework for analysis with infinitesimals.

\subsection{Are there gaps in Euclid?}

Some of the best work on ancient Greece would possibly fail to satisfy
Fried's criteria for authentic history, such as de Risi's monumental
work (\cite{De16}, 2016), devoted to the reception of Euclid in the
early modern age.  Here de Risi writes:
\begin{quote}
Euclid's system of principles has been repeatedly discussed and
challenged: A few gaps in the proofs were found \ldots{}
\cite[p.\;592]{De16}.
\end{quote}
This is a statement about Euclid and not merely its early modern
reception.  Now wouldn't the claim of the existence of a what is seen
today as a ``gap" in Euclid be at odds with Fried's \emph{tabula rasa}
axiom (see Section~\ref{s3})?

\subsection{Philological thought experiment}

Fried's discussion is so general as to raise questions about its
utility.  Dipert notes in his review of the original 1981 edition of
Mueller \cite{Mu06}:
\begin{quote}
It will be difficult in the coming years for anyone doing serious
research on Euclid, \emph{outside of the narrowest philological
studies}, not first to have come to grips with the present book, and
it is to be hoped that this volume will inject new vigor into
discussions of Euclid by contemporary logicians and philosophers of
mathematics.  (Dipert \cite{Di81}, 1981; emphasis on ``philological
studies'' added)
\end{quote}
Inspired by Dipert's observation, we propose the following thought
experiment.  Consider a hypothetical study of, say, the frequency of
Greek roots in the texts of ancient Greek mathematicians.  Surely this
is a legitimate study in Philology.  As far as Fried's requirements
for authentic history, such a study would meet them with flying
colors.  Thus, the satisfaction of the \emph{discontinuity} axiom (see
Section~\ref{s31}) is obvious.  The satisfaction of the \emph{tabula
rasa} axiom (see Section~\ref{s3}) is evident, seeing that no modern
mathematics is used at all in such a study.  The risk of a Platonist
deviation (see Section~\ref{s24}) is infinitesimal.  Freudenthal, van
der Waerden, and Weil may well have written on interpreting the
classics; but by Fried's ideological criteria, our hypothetical
philological study would constitute legitimate mathematical history,
surpassing anything that such ``privileged observers'' may have
written.  Yet it seems safe to surmise, following Dipert, that the
audience for such a philological study among those interested in the
history of mathematics would be limited.

\section{Evolution of Unguru polarity}

Fried admits in his 2018 essay that when he was a graduate student
under Unguru, it was axiomatic that there are only two approaches to
the history of mathemathics: that of a historian, and that of a
mathematician, as we show in Section~\ref{s51}.

\subsection{Fried's upgrade}
\label{s51}

In his 2018 article, Fried writes:
\begin{quote}
By the time I finished my Ph.D., I could make some distinctions: I
could divide historians of mathematics into a mathematician type, such
as Zeuthen or van der Waerden, a historian type, like Sabetai Unguru,
and, perhaps, a postmodern type {\ldots} (\cite[p.\;4]{Fr18})
\end{quote}
The latter ``type'' is quickly dismissed as ``not in fact a serious
option'' leaving us with onlys two options, historian and
mathematician.  Fried goes on to relate in his essay that he came to
appreciate that the historiographic picture is more complex, resulting
in the novel labels of \emph{mathematical conquerors},
\emph{privileged observers}, and the like.  Such a more complex
picture is something of a departure from Unguruan orthodoxy, as we
analyze in Section~\ref{s52}.

\subsection{Unguru polarity}
\label{s52}

In Section~\ref{s51} we summarized Fried's upgrade of Unguru's
original framework.  Meanwhile, Unguru himself sticks to his guns as
far as the original dichotomy of mathematician versus historian is
concerned.  In his 2018 piece, Unguru reaffirms the idea that there
are only two approaches to the history of mathematics:
\begin{quote}
The paper deals with two \emph{polar-opposite} approaches to the study
of the history of mathematics, that of the mathematician, tackling the
history of his discipline, and that of the historian.  (Unguru 2018,
\cite{Un18}, p.\;17; emphasis added)
\end{quote}
Unguru proceeds to reveal further details on the alleged polarity:
\begin{quote}
[S]ince it is always possible to present past mathematics in modern
garb, ancient mathematical accomplishments can be easily made to look
modern and, therewith seamlessly integrated into the growing body of
mathematical knowledge. That this is a historical \emph{calamity} is
not the mathematician's worry
%
%
\ldots\;Never mind that this procedure is tantamount to the
\emph{obliteration} of the history of mathematics as a
\emph{historical} discipline.  Why, after all, should this concern the
mathematician?  \cite[p.\;20]{Un18} (emphasis on \emph{historical} in
the original; emphasis on \emph{calamity} and \emph{obliteration}
added)
\end{quote}
It seems to us that Unguru's assumption that a mathematician does not
care about a possible ``obliteration'' of the history of mathematics
as a historical discipline, is unwarranted.  We note that ``calamity''
and ``obliteration'' are strong terms to describe the work of
respected scholars such as van der Waerden, Freudenthal, and others.
We will examine the issue in more detail in Section~\ref{s43b}.

\subsection{Polarity-driven historiography}
\label{s43b}

What is the driving force behind Unguru's historiographic ideology,
including his readiness to describe the two approaches as
``antagonistic'' (see Section~\ref{s12})?  The ideological polarity
postulated in Unguru's approach appears to involve a perception of
class struggle, as it were, between historians (H-type, our notation)
and mathematicians (M-type, our notation) with their ``antagonistic''
class interests.  As we already noted in Section~\ref{s13}, M-type as
a class does not fare very well relative to the attribute of
\emph{open-mindedness} in the FU ideology.  For a detailed study of
polarity-driven historiography as applied to, or more precisely
against, Felix Klein and (in Unguru's words) ``the obliteration of the
history of mathematics as a historical discipline'' by Mehrtens, see
Bair et al.\;(\cite{18b}, 2017).%
\footnote{Mehrtens (\cite{Me90}, 1990), in an avowedly marxist
approach, postulates the existence of two polar-opposite attitudes
among German mathematicians at the beginning of the 20th century:
modern (M-type, our notation) and countermodern (C-type, our
notation).  Felix Klein had the bad luck of being pigeonholed as a
C-type, along with unsavory types like Ludwig Bieberbach and the
\emph{SS-Brigadef\"uhrer} Theodor Vahlen.  The value of such crude
interpretive frameworks is limited.}

Unguru seeks to forefront the struggle between H-type and M-type as
\emph{the} fundamental ``antagonism'' in terms of which all historical
scholarship must be evaluated, in an attitude reminiscent of the
classic adage ``The history of all hitherto existing society is the
history of class struggles.''  How fruitful is such a historiographic
attitude?  We will examine the issue in the context of a case study in
Section~\ref{s7}.

\subsection{Is exponential notation faithful to Euclid?}
\label{s7}

As a case study illustrating his historiographic ideology, Unguru
proposes an examination of Euclid's Proposition IX.8, dealing with
what would be called today a geometric progression of lengths starting
with the unity.  Unguru already focused on this example over forty
years ago in \cite{Un75}.  Paraphrased in modern terms, the
proposition asserts that in the geometric progression, every other
term is a square, every third term is a cube, etc.  Unguru \cite{Un18}
objects to reformulating the proposition in terms of the algebraic
properties of the exponential notation~$1,a,a^2, a^3,\ldots$, echoing
the criticisms he already made in \cite{Un75}.

Why does Unguru feel that exponential notation must not be used to
reformulate Proposition IX.8?  He provides a detailed explanation in
the following terms:
\begin{quote}
A proposition for the proof of which Euclid has to toil subtly and
painstakingly, and in the course of whose proof he had to rely on many
previous propositions and definitions (e.g., VIII.22 and 23,
def.\;VII.20) becomes a \emph{trivial commonplace}, which is an
immediate outgrowth, a trite after-effect, of our symbolic
notation:~$1$,~$a$,~$a^2$,~$a^3$,
$a^4, a^5, a^6, a^7, \ldots$ As a matter of fact, if we use modern
symbolism, this ceases altogether to be a proposition and its
truthfulness is an immediate and \emph{trivial application} of the
definition of a geometric progression in the particular case when the
first member equals~$1$ and the ratio,~$q=a$, is a positive integer
(for Euclid)!  \cite[p.\;27]{Un18} (emphasis added)
\end{quote}
Unguru claims that using exponential notation causes Euclid's
proposition to become a trivial commonplace severed from Euclid's
``previous propositions,'' and a trivial application of the definition
of a geometric progression.  In this connection, Bl{\aa}sj\"o points
out that
\begin{quote}
Unguru \ldots\;mistakenly believes that certain algebraic insights are
somehow built into the notation itself.  (Bl{\aa}sj\"o \cite{Bl16},
2016, p.\;330)
\end{quote}
Namely,
\begin{quote}
The fact that, for example,~$a^4$ is a square is not by any means
implied by the symbolic notation itself.  The fact
that~$a^{xy}=(a^x)^y$ is a contingent fact, a result that needs
proving.  It is not at all obvious from the very notation itself
{\ldots} (ibid.)
\end{quote}
Thus, contrary to Unguru's claim, Euclid's Proposition IX.8 is
\emph{not} severed from Euclid's ``previous propositions'' which are
similarly more accessible to modern readers when expressed in modern
notation, whose properties require proof just as Euclid's propositions
do.

For instance, Proposition VIII.22 mentioned by Unguru asserts the
following: ``If three numbers are in continued proportion, and the
first is square, then the third is also square.''  In modern
terminology this can be expressed as follows: if~$a^2:b=b:c$, then
$a^2c=b^2$, and therefore~$c=x^2$ for some~$x$.  Put another way,
$a^2rr=(ar)^2$.%
\footnote{Here the first term in the progression is a square~$a^2$ by
hypothesis.  The second term is~$a^2r$ and the third term
is~$(a^2r)r$.  The identity $a^2rr=(ar)^2$ enables one to conclude
that the third term is also a square.}
This is not a triviality but rather an identity that requires proof.
Such an identity could possibly be used in the proof of special cases
of~$a^{xy}=(a^x)^y$.  For more details see Mueller \cite{Mu06}.
Unguru's ideological opposition to using modern exponential notation
in this case has little justification.

It is a pity that (Unguru \cite{Un18}, 2018) chose not to address
Bl{\aa}sj\"o's rebuttal of his objections.  Note that the rebuttal
(Bl{\aa}sj\"o \cite{Bl16}, 2016) appeared two years earlier than
Unguru's piece.

\subsection{What is an acceptable meta-language?}
\label{s61}

In his 2018 piece, Unguru reiterates a sweeping claim he already made
in 1979:
\begin{quote}
The only acceptable meta-language for a historically sympathetic
investigation and comprehension of Greek mathematics seems to be
ordinary language, not algebra.  \cite[p.\;30]{Un18}.
\end{quote}
Given such a stance, it is not surprising that Unguru opposes any and
all use of algebraic notation (including exponential notation) in
dealing with Euclid (see Section~\ref{s7}) and Apollonius (see
Section~\ref{s22}).  However, Berggren notes in his review of (Unguru
\cite{Un79}, 1979) that the reason
\begin{quote}
why modern words, with the concepts they embody, are acceptable as
analytic tools where Renaissance (or even Arabic) algebra is
forbidden,%
\footnote{We added the comma for clarity.}
is never explained [by Unguru].  (Berggren \cite{Be79}, 1979)
\end{quote}
The position of some other historians with regard to Unguru's claims
is discussed in Section~\ref{s7b}.

\section{Do historians endorse Unguru polarity?}
\label{s7b}

Unguru's positing of a polarity of historian \emph{vs} mathematician
tends to obscure the fact that a number of distinguished
\emph{historians} have broken ranks with Unguru on the methodological
issues in question, such as the following scholars.
\begin{enumerate}
\item
Kirsti (M\o ller Pedersen) Andersen wrote a negative review of
Unguru's polarity manifesto \cite{Un75} for Mathematical Reviews,
noting in particular that Unguru ``underestimates the historians'
[e.g., Zeuthen's] understanding of Greek mathematics'' (\cite{An75},
1975).
\item
C. M. Taisbak notes that Unguru and Rowe ``are being ridiculously
unfair, to say the least, towards Heath at this point [concerning
interpretation of the \emph{Elements}, items I44 and I45], to say
nothing of others'' (\cite{Ta81}, 1981).
\item
\'Arp\'ad Szab\'o receives the strongest endorsement from Unguru in
\cite[pp.\;78, 81]{Un75}.  Yet when Szab\'o analyzed \emph{Elements}
Book\;V \cite[p.\;47]{Sz78}, he employed symbolic notation introduced
in the 19th century by Hermann Hankel.%
\footnote{In more detail, Hankel (\cite{Ha74}, 1874, pp.\;389--404)
introduced algebraic notation in an account of Euclid's
\emph{Elements} book\;V.  Furthermore, Heiberg (\cite{He83}, 1883,
vol.\;II, s.\;3) employed Hankel's notation in his translation of
Book\;V into Latin.  For more details see B{\l}aszczyk \cite[p.\;3,
notes 5, 11]{Bl13}.}
Such a practice is clearly contrary to Unguru's position on modern
algebraic notation; see Section~\ref{s61}.  Unlike Unguru, Szab\'o
treats van der Waerden's scholarship with respect and even relies on
it (see Section~\ref{s24}).
\item
Christianidis \cite[p.\;36]{Ch18} proposes a distinction between
premodern algebra and modern algebra and argues that Diophantus can be
legitimately analyzed in terms of the former category (see
Section~\ref{s26}).
\end{enumerate} 

The present article is \emph{not} a defense of the mathematician as
mathematical historian.  The main thrust of this article is the
following.  The postulation of an ideological polarity of historian
\emph{vs} mathematician (the latter routinely suspected of a Platonist
deviation) does more harm than good in that it obscures the only
possible basis for evaluating work in the history of mathematics,
namely competent scholarship.  A mathematician who wishes to write
about a historical figure, but is insufficiently familiar with the
historical period and/or the primary documents, should be criticized
as much as a historian insufficiently familiar with the mathematics to
appreciate the fine points, and indeed the implicit aspects (as
detailed e.g., in Bl{\aa}sj\"o--Hogendijk \cite{Bl18}), of what the
historical figure actually wrote.  The axioms of discontinuity and
tabula rasa and the positing of a polarity between mathematicians and
historians are of questionable value to the task of the history of
mathematics.

\section*{Acknowledgments}  

We are grateful to John T. Baldwin, Viktor Bl{\aa}sj\"o, Piotr
B{\l}aszczyk, Robert Ely, Jens Erik Fenstad, Yvon Gauthier, Karel
Hrbacek, Vladimir Kanovei, and David Sherry for helpful comments on
earlier versions of the article.  We thank Professor A. T. Fomenko,
academician, Moscow State University for granting permission to
reproduce his illustration in Figure~\ref{f1} in Section~\ref{s24}.
The influence of Hilton Kramer (1928--2012) is obvious.

\end{document}